\documentclass[envcountsame,envcountsect]{svmult}
\usepackage{type1cm}
\usepackage{amssymb,bbm,mathtools}
\usepackage[bottom]{footmisc}

\spdefaulttheorem{hypothesis}{Hypothesis}{\bfseries}{\itshape}

\newcommand{\ie}{i.e.}

\newcommand{\ov}{\overline}
\renewcommand{\Re}{\operatorname{Re}}

\newcommand{\supp}{\operatorname{supp}}
\newcommand{\e}{\operatorname{e}}
\newcommand{\iu}{\operatorname{i}}

\newcommand{\ran}{\operatorname{ran}}
\newcommand{\tr}{\operatorname{tr}}

\newcommand{\BVl}{\operatorname{BV}_{\rm loc}}

\newcommand{\sm}[1]{\big(\begin{smallmatrix}#1\end{smallmatrix}\big)}
\newcommand{\bb}[1]{{\mathbb{#1}}}
\newcommand{\mc}[1]{{\mathcal{#1}}}

\newcommand{\id}{\mathbbm 1}

\newcommand{\ol}{{\overline{\lambda}}}
\newcommand{\la}{\lambda}

\newcommand{\mx}{{\rm max}}
\newcommand{\mn}{{\rm min}}

\newcommand{\<}{\langle}
\renewcommand{\>}{\rangle}

\begin{document}
\title*{On the spectral theory of systems of first order equations with periodic distributional coefficients}
\titlerunning{Periodic distributional coefficients}
\author{Kevin Campbell and Rudi Weikard}
\institute{Kevin Campbell \at Department of Mathematics, University of Alabama at Birmingham, Birmingham, AL 35294, USA, \email{campke@uab.edu}
\and
Rudi Weikard \at Department of Mathematics, University of Alabama at Birmingham, Birmingham, AL 35294, USA, \email{weikard@uab.edu}}
\motto{In memory of Sergey Nikolaevich Naboko (1950--2020)}

\maketitle

\abstract{We establish a Floquet theorem for a first-order system of differential equations $u'=ru$ where $r$ is an $n\times n$-matrix whose entries are periodic distributions of order $0$.
Then we investigate, when $n=1$ and $n=2$, the spectral theory for the equation $Ju'+qu=wf$ on $\bb R$ when $J$ is a real, constant, invertible, skew-symmetric matrix and $q$ and $w$ are periodic matrices whose entries are real distributions of order $0$ with $q$ symmetric and $w$ non-negative.}

\section{Introduction}
Periodic structures and periodic phenomena have always played a large role in the sciences and in mathematics.
In 1883 Floquet \cite{ASENS_1883_2_12__47_0} gave a canonical form of the solutions of an $m$-th order homogeneous differential equation with periodic coefficients.
Later his result turned out to be instrumental in the understanding of the associated spectral theory.
Such results are now classical even in the somewhat more general case of a first order system $u'=Au$ with a periodic locally integrable matrix $A$.
There are many excellent sources for these matters but we have benefited most from the books by Eastham \cite{MR3075381} and Brown, Eastham, Schmidt \cite{MR2978285}.
Both have extensive lists of references to further literature on the subject.

In this paper we generalize some of these classical results by allowing the coefficients of the differential equation to be periodic distributions of order $0$.\footnote{Recall that distributions of order $0$ are distributional derivatives of functions of locally bounded variation and hence may be thought of, on compact subintervals of $\bb R$, as measures. For simplicity we might use the word measure instead of distribution of order $0$ below.}
In Section~2 we recall the concept of periodicity for distributions and state its most important properties.
This includes the relationship of distributions of order $0$ with measures.
In Section 3 we state and prove the generalization of Floquet's theorem.
In Section 4 we recall some basic facts of the spectral theory for the case of distributional coefficients.
These are taken from Ghatasheh and Weikard \cite{MR4047968} where one may also find additional information on the history and background of the subject.
In Sections 5 and 6 we specialize to the cases of first and second order systems, respectively.
Some simple examples are briefly considered in Section 7.

We end this introduction with a few words on notation.
The set of complex-valued functions of locally bounded variation on $\bb R$ is represented by $\BVl(\bb R)$.
Any $f\in\BVl(\bb R)$ has right- and left-hand limits denoted by $f^\pm$, respectively.
We also use $f^\#$ for the function $(f^++f^-)/2$ which we call \emph{balanced}.
Corresponding to these kinds of functions we have the subspaces $\BVl^\pm(\bb R)$ and $\BVl^\#(\bb R)$ of $\BVl(\bb R)$.
Identity operators are denoted by $\id$ and $\chi_E$ is the characteristic function associated with the set $E$.
A function $R\in\BVl(\bb R)$ generates a complex measure $dR$ at least on compact subsets of $\bb R$.
The corresponding total variation measure is denoted by $|dR|$.
In particular, for Lebesgue measure we use the symbol $dx$, regarding the symbol $x$ as the identity function.
We will often write $\int f$ in place of $\int f dx$, \ie, when an integral does not explicitly specify a measure it may be taken for granted that integration is with respect to Lebesgue measure.
Similarly, if a range for the integration is not specified, integration is over the whole real axis.

\section{Basic properties of periodic distributions}
A distribution $r$ is a linear functional on the set of test functions, \ie, the set of compactly supported infinitely often differentiable functions from $\bb R$ to $\bb C$, satisfying the following property: for every compact subset $K$ of $\bb R$ there are numbers $C\geq0$ and $k\in\bb N_0$ such that
\begin{equation}\label{dist}
|r(\phi)|\leq C \sum_{j=0}^k \|\phi^{(j)}\|_\infty
\end{equation}
whenever $\phi$ is a test function whose support is contained in $K$.

For example, if $x_0\in\bb R$ is fixed, we have the Dirac distribution defined by $\delta_{x_0}:\phi\mapsto\phi(x_0)$.
Also, for any $f\in L^1_{\rm loc}(dx)$ a distribution $\mathfrak{f}$ is defined by $\phi\mapsto\int f(x)\phi(x)dx$.
Below we will often follow the ubiquitous convention to use the same symbol for $\mathfrak{f}$ and $f$.
Context will serve to distinguish the two meanings.

If we define $r'(\phi)=-r(\phi')$ it follows that $r'$ is again a distribution, called the derivative of $r$.
Distributions also have antiderivatives and we have the following important lemma.
\begin{lemma}[Du Bois-Reymond] \label{dBR}
Suppose the derivative of the distribution $r$ is zero.
Then $r$ is the constant distribution, \ie, there is a complex number $C$ such that $r(\phi)=C\int\phi$ for every test function $\phi$.
\end{lemma}

The number $\omega\in\bb R$ is called a \emph{period} of the distribution $r$, if $r(\phi)={r(\phi(\cdot+\omega))}$ for every test function $\phi$.
The distribution $r$ is called \emph{periodic} (or, more specifically, $\omega$-periodic) if it has a non-zero period $\omega$.
Of course, $0$ is a period of any distribution.

\begin{theorem}\label{t2.2}
A periodic distribution $r$ has the following properties:
\begin{enumerate}
\item If $\omega_1$ and $\omega_2$ are periods of $r$, then so are $\omega_1\pm\omega_2$.
\item If $\omega$ is a period of $r$, then any integer multiple of $\omega$ is also a period.
\item If $r$ is constant, then all real numbers are periods of $r$.
\item If the set of periods of $r$ has a finite limit point, then $r$ is constant.
\item The infimum $\omega_0$ of the set of positive periods of $r$ is itself a period.
If $\omega_0>0$, every period of $r$ is an integer multiple of $\omega_0$.
The period $\omega_0$ is then called the fundamental period of $r$.
\item If $R$ is an antiderivative of $r$, then there is a complex number $\alpha$ and a periodic distribution $P$ with the same periods as $r$, such that $R(\phi)=\alpha\int x\phi(x)\, dx+P(\phi)$.
\end{enumerate}
\end{theorem}

\begin{proof}
Properties (1) -- (3) are trivial.

To prove (4) notice first that $0$ must also be a limit point of the periods of $r$.
We show that $r'=0$, since our claim follows then from du Bois-Reymond's lemma.
Suppose $\omega$ is a period of $r$ with $0<\omega<1$ and let $\phi$ be a test function.
Pick a compact interval $K$ such that the supports of both $\phi$ and $\phi(\cdot+\omega)$ are in $K$.
Then
$$0=\frac{r(\phi(\cdot+\omega))-r(\phi)}{\omega}=r\bigg(\frac{\phi(\cdot+\omega)-\phi}{\omega}\bigg)$$
and hence, for appropriate numbers $C$ and $k$,
$$|r'(\phi)|=|r(\phi')|=|r(\psi)|\leq C \sum_{j=0}^k \|\psi^{(j)}\|_\infty$$
where
$$\psi=\frac{\phi(\cdot+\omega)-\phi}{\omega}-\phi'.$$
Using the mean value theorem twice, we obtain
$$\psi^{(j)}(x)=\phi^{(j+1)}(c)-\phi^{(j+1)}(x)=(c-x)\phi^{(j+2)}(\tilde c)$$
for some $c\in(x,x+\omega)$ and some $\tilde c\in(x,c)$.
Hence $|r'(\phi)|\leq C(k+1)M \omega$ where $M=\max\{\|\phi^{(j+2)}\|_\infty:0\leq j\leq k\}$.
Since $\omega$ may be arbitrarily small, we find $r'(\phi)=0$ and since $\phi$ was arbitrary, we get $r'=0$ as promised.

Property (5) is clear when the infimum $\omega_0$ is $0$, so assume it is not.
Then $r$ is not constant and the previous result shows that $\omega_0$ is not a limit point of the set of positive periods. Instead it must be a period itself.
Now suppose $\omega$ is any other period of $r$.
Then $\omega=n\omega_0+b$ for some $n\in\bb Z$ and $b\in[0,\omega_0)$.
Since $b$ is also a period it must be equal to $0$.

Finally, for property (6) we first obtain from du Bois-Reymond's lemma and the periodicity of $r$ that $R(\phi(\cdot+\omega))-R(\phi)=C\int\phi$ for some number $C$.
Now define the distribution $P$ by setting
$$P(\phi)=R(\phi)+\frac{C}{\omega} \int x\phi(x)dx.$$
Since $\int x\phi(x+\omega) dx=\int(x-\omega)\phi(x) dx$
we find that $P$ is periodic.
Hence the claim follows if we set $\alpha=-C/\omega$.
\end{proof}

\begin{remark}
A function $f\in L^1_{\rm loc}(dx)$ is called periodic with period $\omega$, if ${f(x+\omega)}=f(x)$ for almost all (with respect to Lebesgue measure) $x\in\bb R$.
As mentioned earlier, such a function gives rise to a distribution $\mathfrak{f}$ by setting $\mathfrak{f}(\phi)=\int f\phi$ for any test function $\phi$.
Since
$$\mathfrak{f}(\phi(\cdot+\omega))=\int f\phi(\cdot+\omega)=\int f(\cdot-\omega)\phi=\int f\phi=\mathfrak{f}(\phi)$$
we see that $\mathfrak{f}$ is a periodic distribution with the same periods as $f$.
\end{remark}

In the following we will be concerned only with distributions of order $0$, \ie, those for which one may choose $k=0$ in inequality \eqref{dist} regardless of $K$.
They are in close correspondence with functions of locally bounded variation.
Specifically, if $R\in\BVl(\bb R)$, then it generates a (Borel) measure $dR$ on compact subsets of $\bb R$.
It follows that $\phi\mapsto \int \phi dR$ is a distribution of order $0$, in fact, it is the derivative of the distribution $\phi\mapsto \int R \phi dx$.
Conversely, if $r$ is a distribution of order~$0$, then Riesz's representation theorem shows that there is a function $R\in\BVl(\bb R)$ yielding $r(\phi)=\int \phi dR$.
For brevity we will frequently identify the distribution $r$ and the local measure $dR$.
We will also identify the antiderivative of $r$ with the corresponding function $R$ in $\BVl(\bb R)$.
In particular, we use the designations $r(\phi)$, $\int r\phi$, and $\int \phi dR$ interchangeably.

If $f\in L^1_{\rm loc}(|dR|)$ and $r$ is a distribution of order $0$ we may define the product of $r$ and $f$ (or $f$ and $r$) by setting
$$\phi\mapsto(rf)(\phi)=(fr)(\phi)=\int f\phi dR=\int rf\phi.$$
$rf$ is again a distribution of order $0$.

We need the following substitution rule when dealing with integrals.
\begin{lemma}\label{l3.1}
Suppose $(a,b)$ and $(\alpha,\beta)$ are real intervals and $R:(\alpha,\beta)\to\bb C$ is left-continuous and of bounded variation.
If $T:(a,b)\to(\alpha,\beta)$ is continuous and bijective (and hence strictly monotone), then $R\circ T:(a,b)\to\bb C$ is also left-continuous and of bounded variation.
Moreover, if $g\in L^1(|dR|)$, then
$$\int gdR=\pm \int g\circ T\, d(R\circ T)$$
where one has to choose the positive sign if $T$ is strictly increasing and the negative sign if it is strictly decreasing.
\end{lemma}

This lemma has the following consequence in the context of periodic distributions.
\begin{theorem}\label{t2.5}
Suppose $w$ is a periodic distribution of order $0$ with period $\omega$ and $f\in L^1(|w|)$.
Then $\int wf=\int wf(\cdot+\omega)$.
\end{theorem}
\begin{proof}
Let $T:\bb R\to \bb R: x\mapsto x+\omega$ and let $W$ be an anti-derivative of $w$.
Note that by property (6) of Theorem \ref{t2.2} we have $W(x)=\alpha x+P(x)$ for some periodic function (of locally bounded variation) $P$.
Then $W(T(x))=W(x)+\alpha \omega$ and hence $dW=\alpha+dP=d(W\circ T)$.
\end{proof}

\section{Floquet theory}
In this section we shall develop a Floquet theory for the differential equation
$$u'=ru$$
where $r$ is an $n\times n$-matrix whose entries are periodic distributions of order $0$ all of which have a common period $\omega$
(in this case we call $r$ periodic with period $\omega$).
We seek solutions among \emph{balanced} $\bb C^n$-valued functions of locally bounded variation.

\begin{theorem}\label{t3.2}
Suppose $u$ is a balanced function of locally bounded variation such that $u'=ru$.
If $v$ is defined by $v(x)=u(x+\omega)$, then we also have $v'=rv$.
\end{theorem}

\begin{proof}
Let $T(x)=x-\omega$, let $\phi$ be a test function, and set $\psi=\phi\circ T=\phi(\cdot-\omega)$.
This and Theorem \ref{t2.5} give
$$(rv)(\phi)=\int rv\phi =\int r(v\phi)\circ T = \int r u\psi=(ru)(\psi).$$
We also have, by the translation invariance of Lebesgue measure,
$$v'(\phi)=-\int v(x)\phi'(x) dx=-\int u(x)\psi'(x) dx=u'(\psi).$$
Since the rightmost expressions are the same so are the leftmost.
\end{proof}

Thus the operator which assigns $u(\cdot+\omega)$ to $u$ is a map from the space of solutions of $u'=ru$ to itself.
It is called the \emph{monodromy operator}.

The examples $r=\sm{2&0\\ 0&-2}\sum_{k\in\bb Z} \delta_k$ and $r=\sm{2&0\\ 0&-2}\sum_{k\in\bb Z} (\delta_{2k+1}-\delta_{2k})$, which are periodic distribution with period $1$ and $2$, respectively, show that the solution space of $u'=ru$ may not be $n$-dimensional.
Indeed, in the former case the solution space is trivial while, in the latter case it is infinite-dimensional.
This is due to the fact, that the existence and uniqueness theorem for initial value problems fails for these equations.

To proceed we define the matrix $\Delta_r(x)=R^+(x)-R^-(x)$ and add the following hypothesis.
\begin{hypothesis}\label{h3.2}
Let $\omega>0$.
Assume $r$ is an $n\times n$-matrix of $\omega$-periodic distributions of order $0$ such that the matrices ${\id\pm\frac12 \Delta_r(x)}$ are invertible for every $x\in\bb R$.
\end{hypothesis}

It was shown in \cite{MR4047968} that, under this hypothesis, existence and uniqueness of balanced solutions of initial value problems holds.
It follows immediately that the solution space of $u'=ru$ is $n$-dimensional and hence that we have a fundamental matrix $U$ of solutions.
The determinant of $U(x)$ is different from $0$ for any $x\in\bb R$.
Theorem \ref{t3.2} shows that $U(\cdot+\omega)$ is also a fundamental matrix of solutions.
Hence there is a constant matrix $M$ such that
$$U(x+\omega)=U(x)M.$$
The matrix $M$, called a \emph{monodromy matrix}, depends on the choice of $U$.
If $V$ is another fundamental matrix of solutions so that $V=US$ for a constant invertible matrix $S$ and $\tilde M$ is the associated monodromy matrix, then $M=S\tilde M S^{-1}$, \ie, $M$ and $\tilde M$ are similar matrices.
In particular, they have the same eigenvalues.

It is known from Linear Algebra that we may choose $S$ so that
$\tilde M$ is a matrix in Jordan normal form, \ie, $\tilde M$ is a block diagonal matrix where the diagonal blocks, called Jordan blocks, are square matrices of the form
$$\tilde M_k=\begin{pmatrix}
\rho_k&1&0&\cdots &0\\
0&\rho_k&1&\cdots &0\\
\vdots& \vdots&\ddots &\vdots&\vdots\\
0&0&\cdots&\rho_k&1\\
0&0&\cdots&0&\rho_k\\
\end{pmatrix}.$$
Here $\rho_k$ is an eigenvalue of $\tilde M$ (and of $M$) and hence non-zero.
If the number of Jordan blocks is $s$ and if their size is $\mu_k$, $k=1,...,s$, we have, of course $\sum_{k=1}^s \mu_k=n$.

\begin{theorem}\label{t3.4}
Suppose $r$ satisfies Hypothesis \ref{h3.2}.
The differential equation $u'=ru$ has a fundamental system of solutions of the form
$$\e^{\alpha_k x} \sum_{j=0}^\ell q_{k,j}(x) p_{k,\ell-j}(x) \;\;\text{for $\ell=0, ..., \mu_k-1$ and $k=1,...,s$}$$
where $\e^{\alpha_k \omega}=\rho_k$, $q_{k,0}(x)=1$, $q_{k,j+1}(x)=q_{k,j}(x)\frac{x-j\omega}{(j+1) \rho_k\omega}$, and the $p_{k,j}$ are balanced, periodic $\bb C^n$-valued function with period $\omega$.
\end{theorem}

\begin{proof}
We are adapting Hochstadt's proof in \cite{MR0379932} which avoids introducing the logarithm of $M$.\footnote{There appears to be a flaw in Hochstadt's reasoning which we tried to circumvent.}
Let $v_m$ be the $m$-th column of $V$, the fundamental matrix whose monodromy matrix is in Jordan normal form.
There are unique numbers $k\in\{1,...,s\}$ and $\ell\in\{0,..., \mu_k-1\}$ such that $m=1+\ell+\sum_{h=1}^{k-1}\mu_h$.
Let $m_0=1+\sum_{h=1}^{k-1}\mu_h$.
Then we will prove, by induction over $\ell$, that there are $\omega$-periodic functions $p_{k,j}$ such that
\begin{equation}\label{eq:210921.1}
v_{m_0+\ell}(x)=\e^{\alpha_k x} \sum_{j=0}^\ell q_{k,j}(x) p_{k,\ell-j}(x)
\end{equation}
for $\ell=0, ..., \mu_k-1$.
If $\ell=0$ define $p_{k,0}$ by $p_{k,0}(x)= v_{m_0}(x)\e^{-\alpha_k x}$.
Then the identity $v_{m_0}(x+\omega)=\rho_k v_{m_0}(x)$ shows that $p_{k,0}$ is $\omega$-periodic, \ie, $v_{m_0}$ has the required form.
Now assume that \eqref{eq:210921.1} has been established for $\ell=0, ..., r-1$ including the periodicity of $p_{k,0}$, ..., $p_{k,r-1}$ for some $0<r<\mu_k$.
Let $\ell=r$ and note that $v_{m_0+\ell}$ defines the yet undetermined function $p_{k,\ell}$.
It is only left to prove that $p_{k,\ell}$ is $\omega$-periodic.
Since $v_{m_0+\ell}(x+\omega)-\rho_k v_{m_0+\ell}(x)-v_{m_0+\ell-1}(x)=0$ we get, using the periodicity of $p_{k,0}$, ..., $p_{k,\ell-1}$,
$$p_{k,\ell}(x+\omega)-p_{k,\ell}(x)=\sum_{j=1}^\ell \bigg(\frac1\rho q_{k,j-1}(x)+q_{k,j}(x)-q_{k,j}(x+\omega)\bigg)p_{k,\ell-j}(x).$$
For our choice of the polynomials $q_{k,j}$ induction shows that each term on the right-hand side is $0$ proving the periodicity of $p_{k,\ell}$.
\end{proof}

The eigenvalues $\rho_k$ of a monodromy matrix are called \emph{Floquet multipliers} while the numbers $\alpha_k$ are called \emph{Floquet exponents}.
The associated (generalized) eigenfunctions are called (generalized) \emph{Floquet solutions} of $u'=ru$.

\section{Spectral theory}
The main goal of this paper is to investigate the spectral theory associated with a periodic first order $2\times 2$-system of differential equations.
To this end we will recall some basic definitions and results from \cite{MR4047968} and \cite{MR4298818}.

If $w$ is a non-negative $n\times n$-matrix whose entries are distributions of order $0$, $\tr w$ represents a positive scalar measure.
By $\mc L^2(w)$ we denote the collection of $\bb C^n$-valued functions $f$ whose components are measurable with respect to $\tr w$ and which satisfy $\|f\|^2=\int f^*wf<\infty$.
Then $L^2(w)$ designates the corresponding Hilbert space, \ie, the quotient of $\mc L^2(w)$ by the kernel of $\|\cdot\|$.
The inner product of $L^2(w)$ is, of course, given by $\<f,g\>=\int f^*wg$.
Now consider the differential equation
\begin{equation}\label{de}
Ju'+qu=wf
\end{equation}
where $J$ is a constant, invertible and skew-hermitian $n\times n$-matrix and $q$ is a hermitian $n\times n$-matrix whose entries are distributions of order $0$.
Define the linear relations
$$\mc T_\mx=\{(u,f)\in \mc L^2(w)\times \mc L^2(w): u\in\BVl^\#(\bb R)^n, Ju'+qu=wf\}$$
and
$$\mc T_\mn=\{(u,f)\in \mc T_\mx: \text{$\supp u$ is compact in $\bb R$}\}.$$
Then, in the Hilbert space setting, we represent our differential equation by the relations
$$T_\mx=\{([u],[f])\in L^2(w)\times L^2(w):(u,f)\in\mc T_\mx\}$$
and
$$T_\mn=\{([u],[f])\in L^2(w)\times L^2(w):(u,f)\in\mc T_\mn\}.$$

The cornerstone of spectral theory is the result that $T_\mn^*=T_\mx$, \ie, that $T_\mn$ is a symmetric relation (see \cite{MR4298818}).
As a consequence we have that
$$T_\mx=\ov{T_\mn}\oplus D_{\iu}\oplus D_{-\iu}$$
where $D_\la$ is defined as $\{(u,\la u)\in T_\mx\}$.
These are called \emph{deficiency spaces} if $\la\not\in\bb R$.
The numbers $n_\pm=\dim D_{\pm\iu}$ are called \emph{deficiency indices}.
It is important to recall that $\dim D(\la)$ is independent of $\la$ as long as $\la$ varies in either the upper or the lower half of the complex plane.
The deficiency indices are finite (see \cite{RW20-1}) and if they are identical then there exist self-adjoint restrictions of $T_\mx$ (possibly $T_\mx$ itself).

Even in the case of constant coefficients two complications arise.
Firstly, the space $\mc L_0=\{u: \text{$Ju'+qu=0$ and $wu=0$}\}$ may be non-trivial (note here that $wu=0$ if and only if $\|u\|=0$).
If this happens the problem is called non-definite.
The other issue is that $T_\mx$ may indeed not be a linear operator as our introduction of linear relations already insinuates.

If $T$ is a self-adjoint restriction of $T_\mx$ one defines the resolvent set of $T$ by
$$\varrho(T)=\{\la\in\bb C:  \ker(T-\la)=\{0\}, \ran(T-\la)=L^2(w)\}$$
and the spectrum of $T$ by $\sigma(T)=\varrho(T)^{\rm c}$, the complement of $\varrho(T)$.
For $\la\in\varrho(T)$ the linear relation $(T-\la)^{-1}$ is, in fact, a linear operator from $L^2(w)$ to the domain of $T$.

Define the space $\mc H_\infty=\{f\in L^2(w): (0,f)\in T\}$ and $\mc H_0$ to be its orthogonal complement.
Then the domain of $T$ is a dense subset of $\mc H_0$ and if $(u,f)\in T$, then, of course, $f=f_0+f_\infty$ with $f_0\in \mc H_0$ and $f_\infty\in \mc H_\infty$.
But since $(0,f_\infty)$ is in $T$ so is $(u,f_0)$.
Since $f_0$ is uniquely determined by $u$ we have that $T_0=T\cap(\mc H_0\times\mc H_0)$ is a densely defined self-adjoint linear operator, called the operator part of $T$.
In particular, $\mc H_0=\{0\}$ if and only if the spectrum of $T$ is empty.

\section{The case $n=1$}
\begin{hypothesis}\label{H:4.1}
Throughout this section we assume that $\omega$ is a positive real number and $J$ a non-zero purely imaginary number.
Moreover, $q$ is a real distribution of order $0$ while $w$ is a non-negative but non-zero distribution of order $0$.
Both $q$ and $w$ are periodic with period $\omega$.
\end{hypothesis}

When $\la$ and $x$ are in $\bb R$ the imaginary part of the number
$$B_\pm(x,\la)=J\pm\frac12(\Delta_q(x)-\la\Delta_w(x))$$
is equal to that of $J$ and hence non-zero.
Hypothesis \ref{h3.2} is therefore satisfied for $n=1$.

Let $U(\cdot,\la)$ be the unique solution of the homogeneous equation $J u'+qu=\la wu$ satisfying the initial condition $U^+(0,\la)=1$.
Then $\rho(\la)=U^+(\omega,\la)$ is the Floquet multiplier.
Thus $U(x+n\omega,\la)=\rho(\la)^n U(x,\la)$ whenever $n\in\bb Z$.
Recall from Lemma 3.2 in \cite{MR4047968} that
$$U^+(x,\ol)^*JU^+(x,\la)=U^-(x,\ol)^*JU^-(x,\la)=J.$$
Hence $U^+(x,\ol)^*U^+(x,\la)=1$ and, in particular, $\ov{\rho(\ol)}\rho(\la)=1$.

\begin{theorem}\label{t5.2}
Assume that Hypothesis \ref{H:4.1} holds.
Then the linear relation $T_\mx$ is self-adjoint.
Its spectrum is purely continuous and fills the entire real line.
\end{theorem}

\begin{proof}
Fix $\la$ in $\bb C$ and note that $\rho(\la)\neq0$.
Using Theorem \ref{t2.5} we obtain, for any $n\in\bb Z$,
$$\int_{[n\omega,(n+1)\omega)}|U(\cdot,\la)|^2 w=|\rho(\la)|^{-2n}\int_{[0,\omega)}|U(\cdot,\la)|^2w.$$
Hence
$$\|U(\cdot,\la)\|^2=\sum_{n\in\bb Z} |\rho(\la)|^{-2n} \int_{[0,\omega)} |U(\cdot,\la)|^2w$$
which is finite only when $w=0$, a case we have excluded for being trivial.
Hence no $\la$ can be an eigenvalue of $T_\mx$, \ie, the deficiency spaces are trivial and $T_\mx$ is self-adjoint.

Now assume that $\la$ is real.
Since $\ov{\rho(\ol)}\rho(\la)=1$ this shows that $|\rho(\la)|=1$.
It follows that $\la$ is an element of the so called stability set $S$, the set of those $\la$ for which $U(\cdot,\la)$ is bounded.
We will prove that $S\subset\sigma$, the spectrum of $T_\mx$.
Then
$$\bb R\subset S \subset \sigma \subset \bb R$$
which entails that $S=\sigma=\bb R$.

Thus assume now that $\la\in S$ and, by way of contradiction, that it is also in the resolvent set of $T_\mx$.
If we can construct a sequence $n\mapsto (\phi_n,\la \phi_n+f_n)\in T_\mx$ such that $\|\phi_n\|=1$ and $\lim_{n\to\infty}\|f_n\|=0$, then $\phi_n=R_\la f_n$.
Since $R_\la$, the resolvent operator for $T_\mx$ at $\la$, is bounded\footnote{Even in the case of a relation the resolvent is necessarily an operator.}, we get $1=\|\phi_n\|\leq \|R_\la\| \|f_n\|\to0$, a contradiction.

Let us construct the required sequence in $T_\mx$.
Since $\la$ is fixed we will, in the course of this proof, simply write $U$ in place of $U(\cdot,\la)$.
Let $W$ be the left-continuous anti-derivative of $w$ which vanishes at $0$.
For $n\in\bb Z$ let $I_n=(n\omega,(n+1)\omega]$ and for $n\in\bb N$ define
$$h_n=(W+(n+1)W(\omega))\chi_{I_{-n-1}}+W(\omega)\chi_{(-n\omega,n\omega]}-(W-(n+1)W(\omega))\chi_{I_n}$$
and $\phi_n=a_n (h_nU)^\#$ where the numbers $a_n$ are chosen so that $\|\phi_n\|=1$.
Note that the functions $h_n$ are left-continuous, i.e, $h_n=h_n^-$.

Recall the product rule for functions of locally bounded variation, \ie, $(fg)'=f^+g'+f'g^-=f^-g'+f'g^+$.
Hence $\phi_n'=a_n h_n'U^++a_n h_n U'$.
Since $(h_nU)^\#=h_nU+\frac12(h_n^+-h_n^-)U^+$ we get
\begin{multline*}
J\phi_n'+(q-\la w)\phi_n\\ =a_nh_n (JU'+(q-\la w)U)+a_n\big(Jh_n'+\frac12(q-\la w)(h_n^+-h_n^-)\big)U^+.$$
\end{multline*}
The first term on the right is $0$ while the second is a measure supported on the closure of $I_{-n-1}\cup I_{n}$.
More precisely, we have $h_n'=w\chi_{[-n-1,-n)}-w\chi_{[n,n+1)}$, $h_n^+(x)-h_n^-(x)=\Delta_w(x)$ for $x\in [-n-1,-n)$, and $h_n^+(x)-h_n^-(x)=-\Delta_w(x)$ for $x\in [n,n+1)$.
Now observe that the discrete measures $q\Delta_w$ and $\Delta_qw$ are identical so that we get
$J\phi_n'+(q-\la w)\phi_n = w f_n$
when we define
$$f_n = a_n(J+\frac12(\Delta_q-\la \Delta_w))U^+(\chi_{[-n-1,-n)}-\chi_{[n,n+1)}).$$
It follows that $(\phi_n,\la \phi_n+f_n)$ is an element of $T_\mn\subset T_\mx$.

Next we show that the norming constants $a_n$ tend to $0$.
Indeed, using again Theorem \ref{t2.5} and the fact that $|\rho(\la)|=1$,
$$\|\phi_n\|^2\geq |a_n|^2 W(\omega)^2 \int_{[-n\omega,n\omega)} |U|^2 w =2n |a_n|^2 C$$
where $C=W(\omega)^2\int_{[0,\omega)} |U|^2 w$ does not depend on $n$.
Hence, $|a_n|\leq 1/\sqrt{2nC}$.
It is now also clear that $\|f_n\|$ tends to $0$ so that our proof is finished.
\end{proof}

\section{The case $n=2$, real coefficients}
The following assumptions are in force throughout this section.
\begin{hypothesis}\label{H:5.1}
The $2\times 2$-matrix $J$ is real, skew-hermitian, and invertible.
The entries of the $2\times 2$-matrices $q$ and $w$ are real distributions of order $0$ with $q$ hermitian (symmetric) and $w$ non-negative.
Both $q$ and $w$ are $\omega$-periodic where $\omega$ is a positive real number.
Moreover, the matrices
$$B_\pm(x,\la)=J\pm\frac12(\Delta_q(x)-\la\Delta_w(x))$$
are invertible for all $\la,x\in\bb R$.
\end{hypothesis}

We observe that for all $x\in\bb R$ and all $\la\in\bb C$ we have $\det B_+(x,\la)=\det B_-(x,\la)$.
Let $\Lambda$ be the set of all those $\la\in\bb C$ such that, for some $x\in\bb R$, $\det B_\pm(x,\la)=0$.
Our assumptions guarantee that $\Lambda$ does not intersect $\bb R$ and that there are only finitely many elements of $\Lambda$ in any disk of finite radius.
Hence $\Lambda$ is a discrete set.
It is also symmetric with respect to the real axis.
See \cite{RW20-1} for more details.

Furthermore, our condition on $J$ implies that $J=r\sm{0&-1\\ 1&0}$ for some number $r\in\bb R\setminus\{0\}$.
The motivation to restrict our attention to the real case only, is that $v=\ov u$ satisfies $Jv'+(q-\ol w)v=0$, if $u$ satisfies $Ju'+(q-\la{w})u=0$.
In particular, if $U(\cdot,\la)$, for $\la\not\in\Lambda$, is a fundamental matrix for $Ju'+(q-\la w)u=0$, then $\ov{U(\cdot,\la)}$ is a fundamental matrix for $Ju'+(q-\ol w)u=0$.
It follows that, for $\la\in\bb R$, we may choose the fundamental matrix to be real-valued.
Moreover, $r\det U^\pm(\cdot,\la)=-U^\pm_1(\cdot,\ol)^* J U^\pm_2(\cdot,\la)$ when $U_1$ and $U_2$ denote the first and second column of $U$, respectively.
Lemma 3.2 in \cite{MR4047968} implies therefore that $\det U^+(\cdot,\la)=\det U^-(\cdot,\la)$ is constant.
Thus any monodromy matrix $M(\la)$ has determinant $1$ and this, in turn implies that the Floquet multipliers are reciprocals of each other.
Therefore the Floquet multipliers are uniquely determined (up to transposition) by their sum, the \emph{Floquet discriminant}
$$D(\la)=\tr M(\la).$$
Note that $\tr M(\la)$, just like $\det M(\la)$, is invariant under similarity transforms, \ie, it is independent of the specific fundamental matrix used to find $M(\la)$.
Since we may choose $U(\cdot,\la)$ real for $\la\in\bb R$, it follows that $D(\la)$ is then also real.
Moreover, by Theorem 2.7 of \cite{MR4047968} the entries of $U(x,\cdot)$ are analytic in $\bb C\setminus\Lambda$ for any $x\in\bb R$.
Consequently, the Floquet discriminant $D$ is analytic, too.

Our first major goal is to show the absence of eigenvalues of $T_\mx$.
First we show that no point in $\Lambda$ can be an eigenvalue.

\begin{theorem}\label{t6.4}
If $\la\in\Lambda$, then any solution of $Ju'+qu=\la wu$ is identically equal to $0$.
\end{theorem}

\begin{proof}
Since $\la\in\Lambda$ there is a point $x_0$ for which $B_\pm(x_0,\la)$ are not invertible.
Of course $x_0+\omega$ is then also such a point.
First we show that ${\ran B_+(x_0,\la)}$ and ${\ran B_-(x_0,\la)}$ intersect only trivially.
Let us abbreviate $B_\pm(x_0,\la)$ simply by $B_\pm$.
Neither $B_+$ nor $B_-$ is $0$.
Hence there are two vectors $v_+$ and $v_-$ spanning, respectively, their kernels.
Moreover, $0\neq2Jv_+=B_-v_+$ and $0\neq2Jv_-=B_+v_-$.
Assuming now, by way of contradiction and without loss of generality that $B_-v_+=B_+v_-$ we get $Jv_+=Jv_-$ which contradicts the fact that $J$ is injective and shows that $\ran B_+\cap\ran B_-=\{0\}$.

If $u$ solves the equation $Ju'+qu=\la wu$ we must have $B_+(x_0,\la)u^+(x_0)=B_-(x_0,\la)u^-(x_0)$ and hence, by the above, that $B_\pm(x_0,\la)u^\pm(x_0)=0$.
By the same argument we get $B_\pm(x_0+\omega,\la){u^\pm(x_0+\omega)}=0$.
It follows now that $v=(u \chi_{(x_0,x_0+\omega)})^\#$ is also a solution of $Ju'+qu=\la wu$, in fact a solution of finite norm.
If the norm were positive we would have a complex eigenvalue of $T_\mn$ which is impossible.
Therefore $wv=0$ so that $v$ also solves $Ju'+qu=\la w u$ for all $\la\in\bb C$ including $\la=0$.
However, for $\la=0$ solutions of initial value problems are unique.
Hence $\mc L_0$ is trivial and $u$ is identically equal to $0$.
\end{proof}

\begin{lemma}\label{l6.1}
Fix $\la\in\bb C\setminus\Lambda$.
For a solution $u$ of $Ju'+qu=\la wu$ and $n\in\bb Z$ define
$$I_n(u)=\int_{[n\omega,(n+1)\omega)} u^*wu.$$
If $I_n(u)=0$ for two consecutive non-zero integers, then $\|u\|=0$.
Otherwise $u$ has infinite norm.
If, in the former case, $u$ is not a Floquet solution and not equal to $0$, then all Floquet solutions of $Ju'+qu=\la wu$ also have norm $0$.
\end{lemma}

\begin{proof}
Suppose first that we have two linearly independent Floquet solutions $\psi_1$ and $\psi_2$, the former with multiplier $\rho$ and the latter with multiplier $1/\rho$.
Define $N_j^2=\int_{[0,\omega)}\psi_j^*w\psi_j$ with $N_j\geq0$ and $A=\int_{[0,\omega)}\psi_1^*w\psi_2$.
Note that $|A|\leq N_1N_2$.
If $u=\alpha\psi_1+\beta\psi_2$ we obtain
\begin{multline*}
I_n(u)=|\alpha|^2|\rho|^{2n}N_1^2+2\Re(\ov\alpha\beta(\ov\rho/\rho)^nA)+|\beta|^2|\rho|^{-2n}N_2^2\\
 \geq (|\alpha|N_1 |\rho|^n-|\beta|N_2|\rho|^{-n})^2.
\end{multline*}
If $u=0$ there is nothing to prove and thus we may assume that $\alpha$ and $\beta$ are not both $0$.
If one of $\alpha$ and $\beta$ is $0$, then $u$ is a Floquet solution and either all of the $I_n(u)$ are $0$ or all of them are positive.
Correspondingly, $\|u\|=0$ or $\|u\|=\infty$.
If neither $\alpha$ nor $\beta$ is $0$ let us assume that $I_n(u)=I_{n+1}(u)=0$ for some $n\in\bb Z$.
This implies that either $N_1=N_2=0$ or else $|\rho|=1$.
In the former case all solutions of $Ju'+qu=\la wu$ have norm $0$.
In the latter case we have
$$I_n(u)=2N^2(1+\Re(z\rho^{-2n}))$$
where we put $N=|\alpha|N_1=|\beta|N_2>0$ and $z=\ov\alpha\beta A/N^2$.
Since $|z\rho^{-2n}|\leq 1$ and $I_n(u)=I_{n+1}(u)=0$ we get $z\rho^{-2n}=z\rho^{-2n-2}=-1$ which implies that $\rho^2=1$.
But this means that $\psi_1$ and $\psi_2$ have the same Floquet multiplier.
Hence $u$ itself is a Floquet solution and must have norm $0$.

We now consider the case when there is only one linearly independent Floquet solution.
Then $\rho=\pm1$ and we treat the case when $\rho=1$; the other one is similar.
The general solution of $Ju'+qu=\la wu$, according to Theorem~\ref{t3.4}, is
$$u(x)=\frac{\alpha x}{\omega}\psi(x)+\alpha p(x)+\beta \psi(x),$$
where $\psi$ is a Floquet solution and $p$ is some $\omega$-periodic function.
We confine ourselves to the case $\alpha=1$ and define $\varphi=p+\beta\psi$, which is $\omega$-periodic.
For $n\geq1$ we find $I_n(u)\geq (F_n-G)^2$ where
$F_n^2=\int_{[0,\omega)}\frac{(x+n\omega)^2}{\omega^2} \psi^*w\psi$ and $G^2=\int_{[0,\omega)} \varphi^*w\varphi$.
Now $I_n(u)=I_{n+1}(u)=0$ implies that $F_n=G=F_{n+1}$ and hence $0=\int_{[0,\omega)}(2n+1+2x/\omega) \psi^*w\psi$.
From this we conclude $\psi^*w\psi=0$ so that $\|\psi\|=0$, $F_n=G=0$ and $\|u\|=0$.

Now the only thing left to prove is that a solution $u$ which does not have norm $0$ has infinite norm.
This is already known for Floquet solutions so we assume that $u$ is not a Floquet solution.
For $u$ to have finite norm it is necessary that $I_n(u)$ tends to $0$ as $n$ tends to $\infty$ or $-\infty$.
In the presence of two independent Floquet solutions it is thus necessary that $|\alpha|N_1=|\beta| N_2$ and that $|z|=|\ov\alpha\beta A|/N^2=1$.
In this case we have $I_n(u)={2N^2(1+\cos(t_0+2nt))}$ where we set $z=\e^{it_0}$ and $\rho=\e^{-it}$ for some $t_0,t\in\bb R$.
However, the sequence $n\mapsto I_n(u)$ converges to $0$ only if it is identically equal to $0$.
If only one linearly independent Floquet solution exists and $u$ is not a Floquet solution $I_n(u)$ cannot converge to $0$ unless $F_n=G=0$.
\end{proof}

An immediate consequence of the previous lemma is the next theorem.
\begin{theorem}\label{t6.2}
The linear relation $T_\mx$ has no eigenvalues. In particular, it is self-adjoint.
\end{theorem}

\begin{lemma}\label{L6.4}
Fix $\la\in\bb C\setminus\Lambda$. Let $B$ be the set of the vectors $\int U(\cdot,\ol)^*wf$ when $f$ varies among the functions in $L^2(w)$ with support in $[0,2\omega)$.
If $\mc L_0$ is trivial, then $B=\bb C^2$.
\end{lemma}

\begin{proof}
Assume, by way of contradiction, that $B$ is less than two-dimensional.
Since $B=\{0\}$ would imply $w=0$ it follows that, in fact, $B$ is one-dimensional.
Denote the subspaces of $B$ obtained by restricting $f$ to those functions which are supported only on $I_1=[0,\omega)$ or $I_2=[\omega,2\omega)$ by $B_1$ and $B_2$, respectively.
These are also one-dimensional and hence $B_1=B_2=B$.
Choosing $0\neq\alpha\in B^\perp$ and $f\in L^2(w)$ with support in $I_1$ we get, using Theorem \ref{t2.5} and the periodicity of $w$,
\begin{multline*}
0=\alpha^*\int_{I_2} U(\cdot,\ol)^*w f(\cdot-\omega)=\alpha^*\int_{I_1} U(\cdot+\omega,\ol)^*w f\\ =\alpha^*M(\ol)^*\int_{I_1} U(\cdot,\ol)^*w f.
\end{multline*}
This shows $M(\ol)\alpha$ is perpendicular to $B_1$ and hence a multiple of $\alpha$.
Therefore $u=U(\cdot,\ol)\alpha$ is a Floquet solution of $Ju'+qu=\ol wu$ with norm $0$, \ie, $u\in\mc L_0$ which contradicts our hypothesis.
\end{proof}

If our problem is non-definite $T_\mx$ is a particularly simple relation as our next theorem shows.
\begin{theorem}\label{T6.6}
For the relation $T_\mx$ to be equal to $\{0\}\times L^2(w)$ it is necessary and sufficient that $\mc L_0\neq\{0\}$.
In this situation the spectrum of $T_\mx$ is empty.
\end{theorem}

\begin{proof}
To show necessity assume that $T_\mx=\{0\}\times L^2(w)$ and, by way of contradiction, that $\mc L_0=\{0\}$.
For any $f\in L^2(w)$ with support in $[0,2\omega)$ there is a function $u$ of locally bounded variation such that $Ju'+qu=wf$ and $wu=0$.
In $(-\infty,0)$ we have that $u$ must be equal to a linear combination of (generalized) Floquet solutions which we call $\tilde u$ .
Since then, in the terminology of Lemma \ref{l6.1}, $I_n(\tilde u)=0$ for all negative integers $n$ and since $\mc L_0=\{0\}$ we have, in fact, $\tilde u=0$ and hence $u=0$ on $(-\infty,0)$.
Similarly we can show that $u=0$ on $[2\omega,\infty)$.
For $x>0$ the variation of constants formula (see Lemma 3.3 in \cite{MR4047968}) gives
$$u^-(x)=U^{-}(x,0)J^{-1}\int_{[0,x)}U(\cdot,0)^*wf.$$
In particular, $\int_{[0,2\omega)}U(\cdot,0)^*wf=0$ regardless of $f$, contradicting the findings of Lemma \ref{L6.4}.

To show sufficiency we note first that we may assume $w\neq0$ since otherwise $L^2(w)=\{0\}$.
The conclusion will follow from Theorem 7.3 in \cite{MR4047968} once we establish its hypotheses.
We have already shown that $n_\pm=0$ and it is clear that $0<\dim\mc L_0<2$, the former inequality following from the hypothesis that $\mc L_0$ is not trivial and the latter since $w\neq0$.
It remains to show that $\ker\Delta_w(x)\subset \ker\Delta_w(x)J^{-1}\Delta_q(x)$ for all $x\in\bb R$.
Thus fix $x\in\bb R$ and let $\Delta_q(x)=\sm{a&b\\ b&d}$ and $\Delta_w(x)=\sm{\alpha&\beta\\ \beta&\delta}$.
If the rank of $\Delta_w(x)$ is $0$ or $2$, then the needed inclusion holds trivially.
Hence we assume that $\alpha\delta=\beta^2$ with at least one of those numbers non-zero.
In this case $\det B_\pm(x,\la)$ is a real polynomial in $\la$ of degree at most $1$.
Since $\det B_\pm(x,\cdot)$ must not have a real zero, it follows that the coefficient of $\la$ must be $0$, \ie, $a\delta+\alpha d=2b \beta$.
Under that condition it turns out that $\Delta_w(x)J^{-1}\Delta_q(x)$ is symmetric and hence equal to $-\Delta_q(x)J^{-1}\Delta_w(x)$ which proves the desired inclusion and hence the lemma.

For the last claim we simply note that the operator part of $T_\mx$ is $\{([0],[0])\}$ which has no spectrum.
\end{proof}

We now define the conditional stability set $S$ to be the set of those $\la\in\bb C$ such that the equation $Ju'+qu=\la wu$ has a non-trivial bounded solution.
Note that Theorem \ref{t6.4} shows that $S$ and $\Lambda$ do not intersect.
It follows that the conditional stability set is the set of those $\la\in\bb C\setminus\Lambda$ for which the monodromy matrix $M(\la)$ has an eigenvalue of absolute value $1$.

\begin{lemma}\label{l6.6}
If $\la$ is not in $S\cup\Lambda$, then $M(\la)$ has distinct eigenvalues $\rho_1=\e^{\tilde m \omega}$ and $\rho_2=\e^{-\tilde m \omega}$ where $m=\Re\tilde m>0$.
The associated eigenfunctions $\psi_{1,2}$ may be chosen so that  $(\psi_1^{\pm})^\top J\psi_2^\pm=1$.
If, for a given $f\in L^2(w)$, we define
$$u^-(x)=\psi_2^-(x)\int_{(-\infty,x)}\psi_1^\top wf+\psi_1^-(x)\int_{[x,\infty)}\psi_2^\top wf,$$
then $u=u^\#$ satisfies the differential equation $Ju'+qu=w(\la u+f)$.
\end{lemma}

\begin{proof}
The first two statements follow from Theorem \ref{t3.4} and Lemma 3.2 in~\cite{MR4047968}.

As before we shall write $B_\pm$ in place of $B_\pm(\cdot,\la)$.
Since $\psi_j=(\psi_j^++\psi_j^-)/2$, $B_+\psi_j^+=B_-\psi_j^-$, and $B_++B_-=2J$, our hypothesis on the normalization of $\psi_{1,2}$ shows that
\begin{equation}\label{eq20211011.1}
B_+(\psi_2^+\psi_1^\top -\psi_1^+\psi_2^\top)=\id.
\end{equation}
Since
$$u^+(x)=\psi_2^+(x)\int_{(-\infty,x]}\psi_1^\top wf+\psi_1^+(x)\int_{(x,\infty)}\psi_2^\top wf$$
equation \eqref{eq20211011.1} gives $B_+u^+-B_-u^-=\Delta_wf$ and hence
\begin{equation}\label{eq20211011.2}
Ju^-=B_+u^\#-\frac12 \Delta_wf.
\end{equation}
Now recall that, by the definition of Lebesgue-Stieltjes measures, the derivatives of $u^-$ and $u$ are the same.
Hence the product rule for functions of locally bounded variation gives
$$Ju'=J\psi_2'\int_{(-\infty,\cdot)}\psi_1^\top wf+J\psi_1'\int_{[\cdot,\infty)}\psi_2^\top wf
 +J(\psi_2^+\psi_1^\top -\psi_1^+\psi_2^\top)wf.$$
Using the identity $J\psi_k'=(\la w-q)\psi_k=(\la w-q)B_+^{-1}J\psi_k^-$ and equations \eqref{eq20211011.1} and \eqref{eq20211011.2} we obtain now
\begin{multline*}
Ju'=(\la w-q)B_+^{-1}J u^-+JB_+^{-1}wf\\ =(\la w-q)u^\#+\frac12(q-\la w)B_+^{-1}\Delta_wf+JB_+^{-1}wf.
\end{multline*}
Since
$$\frac12(q-\la w)B_+^{-1}\Delta_wf=\frac12(\Delta_q-\la \Delta_w)B_+^{-1}wf=(\id-J B_+^{-1})wf$$
we finally find that $Ju'+(q-\la w)u^\#=wf$.
\end{proof}

Next we show that the function $u$ just constructed is actually an element of $\mc L^2(w)$.
In the following we use the notation of Lemma \ref{l6.6} and its proof freely.
First note that averaging $u^+$ and $u^-$ gives
\begin{multline}\label{eq20211012.1}
u(x)=\psi_2(x)\int^x\psi_1^\top wf+\psi_1(x)\int_x\psi_2^\top wf\\ +\frac12 (\psi_2^+(x)\psi_1(x)^\top+\psi_1^-(x)\psi_2(x)^\top)\Delta_w(x)f(x)
\end{multline}
where we introduced the notation $\int^x$ and $\int_x$ as abbreviations for $\int_{(-\infty,x)}$ and $\int_{(x,\infty)}$, respectively.

Introduce $\tilde W$, the Radon-Nikodym derivative of $w$ with respect to the positive scalar measure $\tr w$.
Then $\tilde W$ is a real symmetric matrix and $0\leq \tilde W\leq\id$ pointwise almost everywhere with respect to $\tr w$ and we may assume it holds indeed everywhere.
Therefore $f^*\tilde W(x) g$ is a semi-inner product in $\bb C^2$ which we denote by $\<f,g\>_x$.
The corresponding semi-norm is $|f|_x=\sqrt{\<f,f\>_x}$.

Recall that $m=\Re\tilde m>0$.
Since $\psi_1(x)=\e^{\tilde mx}p_1(x)$ we get $|\psi_1(x)|_x=\e^{mx} |p_1(x)|_x$.
Similarly $|\psi_2(x)|_x=\e^{-mx} |p_2(x)|_x$.
Since $p_1$ and $p_2$ are periodic we may find a constant $K$ such that $|\psi_j(x)|_x\leq K \e^{\pm mx}$ where the sign in the exponent is $(-1)^{j+1}$.
The norm of a rank one matrix $ab^\top$ is the product of the norms of $a$ and $b$.
Hence we may estimate the $|\cdot|_x$-norm of the last term in \eqref{eq20211012.1} by $K^2 |\Delta_w(x) f(x)|_x$.
However, the matrices $\Delta_w(x)$ have to be uniformly bounded by, say, a number $D\geq 1$ since $w$ is periodic and locally finite.
Thus our estimate on the last term in \eqref{eq20211012.1} becomes $K^2 D|f(x)|_x$.

To deal with the other terms in \eqref{eq20211012.1} we define
$$\phi_1(x)=\int^x \psi_1^\top wf \quad\text{and}\quad \phi_2(x)=\int_x \psi_2^\top wf.$$
Upon using Cauchy's inequality for the $\<\cdot,\cdot\>_x$-inner product we obtain the estimates $|\phi_j(x)|\leq K H_j(x)$
where
$$H_1(x)=\int^x\e^{my} |f(y)|_y \tr w(y) \quad\text{and}\quad H_2(x)=\int_x \e^{-my} |f(y)|_y \tr w(y).$$
Therefore we find
$$|u(x)|_x \leq K^2\e^{-mx}H_1(x) + K^2\e^{mx}H_2(x) + K^2D |f(x)|_x.$$

Since $\|u\|^2=\int |u(x)|_x^2 \tr w(x)$ and $H_j^2\leq H_j^{2\#}$ we obtain
$$\|u\|^2\leq 3K^4D^2 \int \big(\e^{-2mx}H_1^{2\#}(x)+\e^{2mx}H_2^{2\#}(x)+|f(x)|_x^2\big) \tr w(x).$$
Next we estimate the integrals $I_1(a,b)=\int_{(a,b)} \e^{-2mx}H_1^{2\#}(x) \tr w(x)$ and $I_2(a,b)=\int_{(a,b)} \e^{2mx}H_2^{2\#}(x)\tr{w(x)}$ assuming that $a$ and $b$ are points of continuity of $W$, the antiderivative of $w$.
To this end we define
$$R_-(x)=\int_x \e^{-2my} \tr w(y) \quad\text{and}\quad R_+(x)=\int^x \e^{2my} \tr w(y).$$
We emphasize that $R_\pm$ are positive functions and that $R_+$ is non-decreasing while $R_-$ is non-increasing.
Hence $R_\pm\leq R_\pm^\#$.
Using the periodicity of $w$ and Theorem \ref{t2.5} we find $R_+(x)\leq C_m \e^{2mx}$ and $R_-(x)\leq C_{-m} \e^{-2mx}$ for appropriate constants $C_{\pm m}$.
The continuity of these upper bounds implies that even $R_+^\#(x)\leq C_m \e^{2mx}$ and $R_-^\#(x)\leq C_{-m} \e^{-2mx}$.
Now we see that, using integration by parts,
$$I_1(a,b)=-\int_{(a,b)} H_1^{2\#}dR_-=-H_1(b)^2 R_-(b)+H_1(a)^2 R_-(a)+\int_{(a,b)} R_-^\#dH_1^2.$$
The first term on the right is negative and may be omitted.
To deal with the second term we use Cauchy's inequality to show that $H_1(a)^2 \leq R_+(a) \int^a |f|_y^2 \tr w(y)$.
Hence $\lim_{a\to-\infty}H_1(a)^2R_-(a)=0$.
In the last term we use the chain rule $dH_1^2=2H_1^\# dH_1$ (see Vol$'$pert \cite{MR0216338}, Section 13.2 or Vol$'$pert and Hudjaev \cite{MR785938}, Chapter 4, \S6.3).
Thus
\begin{multline*}
\int_{(a,b)} R_-^\#dH_1^2=2\int_{(a,b)} R_-^\# H_1^\# dH_1\\ \leq 2C_{-m}\int_{(a,b)} \e^{-mx}H_1^\#(x) |f(x)|_x \tr w(x).
\end{multline*}
Cauchy's inequality and the fact that $H_j^{\#2}\leq H_j^{2\#}$ give next that
$$I_1(a,b) \leq H_1(a)^2 R_-(a)+2C_{-m} I_1(a,b)^{1/2} \|f\|.$$
Taking now the limit as $a$ tends to $-\infty$ and $b$ to $\infty$ shows that $I_1(-\infty,\infty)\leq 4C_{-m}^2 \|f\|^2$.
Treating $I_2$ in a similar way we finally obtain
$$\|u\|^2 \leq 3K^4D^2(4C_{-m}^2+4C_m^2+1)\|f\|^2.$$
This proves that $(u,\la u+f)\in T_\mx$, \ie, $\la\in\varrho$.
Equation \eqref{eq20211012.1} means that the Green's function for $T_\mx$ is given by
$$G(x,y)=\begin{cases} \psi_2(x)\psi_1(y)^\top &\text{if $y<x$}\\
\frac12(\psi_2(x)\psi_1(x)^\top+\psi_1(x)\psi_2(x)^\top) &\text{if $y=x$}\\
\psi_1(x)\psi_2(y)^\top &\text{if $x<y$.}\end{cases}$$

We have now shown that $(S\cup\Lambda)^{\rm c}\subset\varrho$.
Since $\Lambda$ is also a subset of $\varrho$ we have established the following lemma.
\begin{lemma}\label{l6.7}
The complement $S^{\rm c}$ of the conditional stability set $S$ is contained in $\varrho$, the resolvent set of $T_\mx$.
\end{lemma}

Next we will address the converse statement.
\begin{lemma}\label{l6.8}
If $\mc L_0$ is trivial the conditional stability set $S$ is contained in $\sigma$, the spectrum of $T_\mx$.
\end{lemma}

\begin{proof}
Suppose $\la\in S$ and, by way of contradiction, that it is also in $\varrho$.
Theorem~\ref{t6.4} shows that $\la\not\in\Lambda$.
Hence there is a balanced, bounded Floquet solution $\psi$ of the differential equation $J\psi'+(q-\la w)\psi=0$ with Floquet multiplier $\rho$.
Of course, $\psi=U(\cdot,\la)c_0$ for some vector $c_0$ which is an eigenvector of $M(\la)$ associated with $\rho$.
Lemma 3.2 in \cite{MR4047968} shows that $Ju_0$ is an eigenvector of $M(\ol)^*$ associated with the eigenvalue $1/\rho$.
Lemma~\ref{L6.4} provides the existence of a function $g_0$ compactly supported in $[0,2\omega)$ satisfying
$$\int_{[0,2\omega)} U(\cdot,\ol)^*wg_0=-Jc_0.$$
For $n\in\bb N$ we put $g_n=\rho^n g_0(\cdot-n\omega)$, a function which is supported in ${[n\omega,(n+2)\omega)}$, and use the variation of constants formula to obtain a solution $v_n$ of the equation $Ju'+qu=w(\la u+g_n)$, \ie,
$$v_n^-(x)=U^-(x,\la)\big(c_0+J^{-1}\int_{[0,x)} U(\cdot,\ol)^*wg_n\big).$$
The function $v_n$ coincides with $\psi$ on $(-\infty,n\omega)$.
Since
\begin{multline*}
\int_{[0,(n+2)\omega)} U(\cdot,\ol)^*wg_n=\rho^n M(\ol)^{*n}\int_{[0,2\omega)} U(\cdot,\ol)^*wg_0\\
 =-\rho^n M(\ol)^{*n}Jc_0=-Jc_0
\end{multline*}
it follows that $v_n$ is identically equal to $0$ beyond $(n+2)\omega$.
With a similar device for the negative half-line we may now construct a sequence ${(u_n,\la u_n+f_n)}\in T_\mx$ with the following properties: (i)  $u_n=0$ outside $[-(n+2)\omega,(n+2)\omega]$, (ii) $u_n=\psi$ on $(-n\omega,n\omega)$, and (iii) $\|f_n\|$ is independent of $n$.

Since $|\rho|=1$ we have that
$$\|u_n\|^2\geq \int_{[-n\omega,n\omega)} \psi^*w\psi =2n C$$
where $C=\int_{[0,\omega)} \psi^*w\psi$ does not depend on $n$.

Arguing similarly to the conclusion in the proof of Theorem \ref{t5.2} we now have $\|u_n\|\leq \|(T-\la)^{-1}\| \|f_n\|$ where the left-hand side tends to infinity while the right-hand side is constant, the desired contradiction.
\end{proof}

Combining Theorem \ref{t6.2}, Lemma \ref{l6.7}, and Lemma \ref{l6.8} we have now the following result.
\begin{theorem}\label{t6.10}
Assume that $\mc L_0$ is trivial.
Then the spectrum $\sigma$ of $T_\mx$ is purely continuous and coincides with the conditional stability set $S$.
\end{theorem}

We end this section by an investigation of the properties of the Floquet discriminant $D$.
Suppose $x_0$ is a point of continuity for $Q$ and $W$, the antiderivatives of $q$ and $w$, and $U(\cdot,\la)$ is a fundamental matrix for $Ju'+qu=\la wu$ such that $U(x_0,\la)=\id$.
Observe that
$$J\dot U'(\cdot,\la)+(q-\la w)\dot U(\cdot,\la)=w U(\cdot,\la)$$
where the accent $\dot{}$ denotes the derivative with respect to $\la$.
Since $M={U(x_0+\omega,\cdot)}$ and $\dot U(x_0,\cdot)=0$ the variation of constants formula shows that
$$\dot M(\la)=U(x_0+\omega,\la)J^{-1}\int_{(x_0,x_0+\omega)} U(\cdot,\ol)^*w U(\cdot,\la)=M(\la)J^{-1}T(\la)$$
where the last equation defines the matrix $T$.

Assume that $\la\in\bb R$.
Then $M(\la)$ and $D(\la)$ are real and $T(\la)$ is real and positive semi-definite.
Since $\det M=1$ one may now check that $D^2-4=(M_{11}-M_{22})^2+4M_{12}M_{21}$.
Moreover, setting $a=(M_{11}-M_{22},2M_{21})^\top$ and $b=(-2M_{12},M_{11}-M_{22})^\top$ one obtains
\begin{equation}\label{20211025.1}
4rM_{21}\dot D=-T_{11}(D^2-4)+a^*Ta \;\;\text{and}\;\; 4rM_{12}\dot D=T_{22}(D^2-4)-b^*Tb.
\end{equation}
Now suppose that $D(\la)^2-4<0$.
It is then clear that $M_{12}(\la)$ and $M_{21}(\la)$ are different from $0$ and thus the above identities show that $\dot D(\la)\neq0$ provided that at least one of $T_{11}(\la)$ and $T_{22}(\la)$ is different from $0$.
But $T_{11}(\la)=T_{22}(\la)=0$ may only happen when $w=0$, a case excluded by our hypothesis.

\begin{lemma}
The Floquet discriminant $D$ is constant if and only if $T_\mx$ equals $\{0\}\times L^2(w)$.
Moreover, if $\dim \mc L_0=1$, $D$ does not take its value in $(-2,2)$.
\end{lemma}

\begin{proof}
Suppose $D$ is constant.
If $D^2>4$, then $\bb R\subset S^c\subset \varrho$, the resolvent set of $T_\mx$.
Since the spectrum is empty the domain of $T_\mx$ must be $\{[0]\}$.
Since $T_\mx$ is self-adjoint its range must be $L^2(w)$.
If $D^2\leq 4$, then $S=\bb C\setminus\Lambda$ in view of the analyticity of $D$.
If $\mc L_0$ were trivial, Lemma \ref{l6.8} would show that $\sigma=\bb C$.
Now Theorem \ref{T6.6} shows that $T_\mx=\{0\}\times L^2(w)$.

Conversely, if $T_\mx=\{0\}\times L^2(w)$ and, consequently, $\mc L_0\neq\{0\}$, there exists a non-trivial solution of $Ju'+qu=0$ with norm $0$.
Lemma \ref{l6.1} shows that then at least one of the Floquet solutions must have norm $0$.
Such a solution is then a solution of $Ju'+qu=\la wu$ for all $\la$.
If $\rho$ is the associated Floquet multiplier, then $D(\la)=\rho+1/\rho$ for all $\la\in\bb C$.

Finally, suppose $\dim \mc L_0=1$ and, by way of contradiction, that $D$ is in $(-2,2)$.
Then at least one of the diagonal entries of the matrix $T$ in \eqref{20211025.1} is positive.
Thus $\dot D(\la)\neq0$ (for $\la\in\bb R$), a contradiction.
\end{proof}

\begin{theorem}
Suppose $\mc L_0$ is trivial.
Then, the Floquet discriminant $D$, restricted to the real line, is strictly monotone in intervals where $D^2<4$.
If $D(\la)^2=4$, we have $\dot D(\la)=0$ if and only if $M(\la)=\pm\id$.
Such a point is not a strict local minimum if $D(\la)=2$ nor a strict local maximum if $D(\la)=-2$.
\end{theorem}

\begin{proof}
We already mentioned above that $\la\mapsto D(\la)$ is real-valued and analytic as $\la$ varies in $\bb R$.
We have also proved that $\dot D(\la)\neq0 $ when $D(\la)^2<4$.

Now, assume $D(\la)^2=4$ and $\dot D(\la)=0$.
Then \eqref{20211025.1} shows that $a^*Ta=b^*Tb=0$.
Since $\mc L_0$ is trivial $T$ is positive definite and hence $a$ and $b$ are $0$ proving $M(\la)=\pm\id$.
Conversely, if $M(\la)=\pm\id$, then $\dot D(\la)=\pm\tr(J^{-1} T(\la))=0$.
Finally, if $D$ has a strict local minimum of $2$ at a point $\la$, then $\la$ is an isolated point in the spectrum of $T_\mx$ and hence an eigenvalue which is impossible. 
Of course, a similar argument prevents $D$ to have a local maximum of $-2$ anywhere.
\end{proof}

\section{Examples}
In all examples presented below we assume that $n=2$, $J=\sm{0&-1\\ 1&0}$ and $\omega=1$ (except in our last example).

Our first three examples have constant coefficients.
\begin{itemize}
\item $q=\sm{0&0\\ 0&-1}$ and $w=\sm{1&0\\ 0&0}$.
Here we have $\mc L_0=\{0\}$ and $D(\la)=2\cos(\sqrt{\la})=2\cosh\sqrt{-\la}$.
Hence $S=\sigma=[0,\infty)$.
This is the system describing the second order equation $-y''=\la y$.

\item $q=\sm{a&b\\ b&d}$ and $w=0$.
We get $D=2\cosh\sqrt{b^2-ad}=2\cos\sqrt{ad-b^2}$.
Hence $D$ is constant (as anticipated) and may take any value in $[-2,\infty)$.
If $ad-b^2>0$ we have $S=\bb R$ but $\sigma=\emptyset$.

\item $q=\sm{a&b\\ b&0}$ and $w=\sm{1&0\\ 0&0}$.
Now $\mc L_0$ is spanned by the vector $(0,\e^{bx})^\top$.
We have $D=2\cosh b$. Note that $D$ cannot be in $(-2,2)$.
\end{itemize}
These last two examples show that the hypothesis of a trivial $\mc L_0$ is necessary in Lemma \ref{l6.8}.

Our next three examples involve discrete measures as coefficients.
Let $\mu=\sum_{k\in\bb Z}\delta_k$ where $\delta_k$ is the Dirac measure concentrated on $\{k\}$.
\begin{itemize}
\item $q=\sm{a \mu&0\\ 0&-1}$ and $w=\sm{\alpha\mu&0\\ 0&0}$.
If $\alpha>0$ this problem is definite, \ie, $\mc L_0$ is trivial.
The discriminant is $D(\la)=2+a-\alpha\la$.
Hence $\sigma=[a/\alpha,(4+a)/\alpha]$.
However, if $\alpha=0$, we have $\dim\mc L_0=2$.
In this case $D$ is constant and the constant may take any real value.

\item $q=\sm{a&b\\ b&0}\mu$ and $w=\sm{1&0\\ 0&0}\mu$ where $b^2\neq4$ to avoid an intersection of $\Lambda$ and $\bb R$.
Here $\mc L_0$ is one dimensional.
We get $D=2(4+b^2)/(4-b^2)$ which may attain any value outside $[-2,2)$.

\item $q=\sm{a&b\\ b&d}\mu$ and $w=\id\mu$ where we require $(a-d)^2+4b^2<16$ to avoid that $\Lambda$ intersects $\bb R$.
Then
$$D(\la)=\frac{16}{\la^2-(a+d)\la+ad-b^2+4}-2.$$
Our condition on $a$, $b$ and $d$ guarantees that the denominator is always positive.
Thus $D(\la)>-2$ for all $\la\in\bb R$ but $D(\la)$ approaches $-2$ as $\la$ tends to $\pm\infty$.
There is one maximum for $D$ at $\la=(a+d)/2$.
Its value is at least $2$ and is equal to $2$ precisely when $a=d$ and $b=0$.
Except for this special situation the spectrum consists of two rays separated by one gap.
\end{itemize}

Finally we mention one example which violates our condition that $\Lambda\cap \bb R=\emptyset$.
The example is taken from \cite{RW20-1} where more details may be found.
\begin{itemize}
\item $q=\sm{0&2\\ 2&0}\sum_{k\in\bb Z} (\delta_{2k}-\delta_{2k+1})$ and $w=\sm{2&0\\ 0&0}\sum_{k\in\bb Z}\delta_k$.
Note that here $\omega=2$.
In this case $\Lambda=\bb C$, \ie, unique continuation of solutions fails for any $\la$.
Since there are compactly supported solutions of $Ju'+qu=0$, it turns out that $0$ is an eigenvalue of infinite multiplicity.
The resolvent set is $\bb C\setminus\{0\}$.
\end{itemize}


\end{document}